\documentclass[11pt]{article}
\usepackage[english]{babel}
 \usepackage{makeidx}
 \usepackage{amssymb,amsbsy,amsmath,amsfonts}
\usepackage{graphics,graphicx}
\usepackage{color}
\def\R{{\mathbb R}}
 
\def\v3{\vskip0.3cm \noindent}

\newtheorem{rem}{\bf Remark}[section]

\newtheorem{thm}{\bf Theorem}

\newtheorem{lem}{\bf Lemma\/}[section]
 \newtheorem{hyp}{\bf Hypothesis\/}

\def\QQ{{\rlap {\raise 0.4ex \hbox{$\scriptscriptstyle  $}}
\hskip -0.2em Q}}
\def\RR{{ I\!\!R}}
\def\R{\RR}

\def\CC{{ \rlap {\raise 0.4ex \hbox{$\scriptscriptstyle  $}}
\hskip -0.2em C}}

\def\v3{\vskip0.3cm \noindent}

\begin{document}
\centerline{\Large Sign of the solution to a non-cooperative system}
\vskip0.3cm \centerline  
{\large B\'en\'edicte Alziary}
\par \centerline
{ \small TSE \& IMT (UMR 5219) - CEREMATH-UT1}
\par \centerline   
{\small  Universit\'e de Toulouse, 31042 TOULOUSE Cedex, France }
\par \centerline   
{\small alziary@ut-capitole.fr} 
\vskip0.3cm  \centerline   
{\large Jacqueline Fleckinger,}
\par \centerline{\small IMT (UMR 5219) - CEREMATH-UT1}
\par \centerline   
{\small Universit\'e de Toulouse, 31042 TOULOUSE Cedex, France}
\par \centerline   
{\small  jfleckinger@gmail.com}


\begin{abstract} 
{\small
Combining  the results of a   recent paper by  Fleckinger-Hern\'andez-deTh\'elin \cite{FHeT2015} for a non cooperative $2\times 2$ system with the method of  PhD Thesis of  MH Lecureux 
we compute the sign of the solutions of a $n\times n$  non-cooperative systems  when the parameter varies near 
the lowest principal  eigenvalue of the system. }

 \end{abstract}


\section{Introduction}
Many results have been obtained since decades on Maximum Principle and Antimaximum principle for second order elliptic partial differential equations involving $e.g.$  Laplacian, p-Laplacian, Schr\" odinger operator, ... or weighted equations. Then most of these results have been extended to systems. 
\par \noindent
The maximum principle (studied since centuries) has many applications in various domains as physic, chemistry, biology,...Usually it shows that for positive data the  solutions are positive (positivity is preserved). It is generally valid for a parameter below the "principal" eigenvalue (the smallest one). 
The Antimaximum principle, introduced in 1979 by Cl\'ement and Peletier (\cite{ClPe}), shows that, for one equation,  as this parameter goes through this principal eigenvalue, the sign are reversed; this holds only for a small interval. The original proof 
relies on a decomposition into the groundstate (principal eigenfunction of the operator) and its orthogonal. It is the same idea which has been used in 
\cite{FHeT2015} (combined with a bootstrap method) to derive a precise estimate for the validity  interval of the Antimaximum principle for one equation. By use of this result,  Fleckinger-Hern\'andez-deTh\'elin (\cite{FHeT2015})  deduce results on the sign  of solution for some $2 \times 2$ non-cooperative systems.  
Indeed many papers have appeared for cooperative systems involving various elliptic operators: (\cite{AFLW}, \cite{AFTa}, \cite{Am2005}, \cite{dFMi1986}, \cite{FiMi1}, 
\cite{FiMi2}, \cite{FGoTaT}, \cite{FHeT}, ...).  Concerning non cooperative systems the literature is more restricted (\cite{CaMi},  \cite{FHeT2015},..). 
\par \noindent
In this paper  we extend the results obtained  in \cite{FHeT2015}, valid for $2 \times 2$ non-cooperative systems  involving Dirichlet Laplacian,  to $n \times n$ ones. 
Recall that a system is said to be "cooperative"  if all the terms outside the diagonal of the associated square matrix are positive. 
\par \noindent
For this aim we combine the precise estimate for the validity interval of the antimaximum principle obtained in \cite{FHeT2015} with the method used in \cite{Le}, \cite{AFLW} for systems. 
\vskip0.2cm \noindent
In Section 2 we are concerned with one equation. We first recall the precise estimate for the validity interval for the antimaximum  principle (\cite{FHeT2015}); then  we  give some related results used in  the study of systems.
\par \noindent
In Section 3 we first state our main results for a $n \times n$  system (eventually non-cooperative) and 
then we prove them. 
\par \noindent
Finally,  in Section 4,  we compare our results with the ones of \cite{FHeT2015}. Our method,  which uses the matricial calculus and in particular Jordan decomposition, allows us to have a more general point of view, even for a $2 \times 2$ system.


\section{ Results for one equation: }$\,$ 
In \cite{FHeT2015}, the authors  consider a non-cooperative $2 \times 2$ system 
with constant coefficients. 
Before studying the system they consider one equation and  establish a  precise estimate 
of the validity interval for the antimaximum principle. We recall this result that we use later. 
\par \noindent
\subsection{A precise Antimaximum for the equation \cite{FHeT2015}} 
Let $\Omega$ be a smooth bounded domain in $I\!\!R^N$. Consider the following Dirichlet boundary value problem
\begin{equation} \label{1}
-\Delta z \, =\, \sigma z + \,h \;\, {\rm in }\; \Omega \; ,  \;\; z=0 \; {\rm on } \; \partial \Omega,  \end{equation}
where $\sigma$ is a real parameter. 
\par \noindent
The associated 
 eigenvalue problem  is 
\begin{equation} \label{2 }-\Delta \phi \, =\, \lambda \phi  \;\, {\rm in }\; \Omega \; , \;\; \phi=0 \; {\rm on } \; \partial \Omega. \end{equation}
As usual, denote by  $0<\lambda_1 < \lambda_2 \leq ...$    the eigenvalues of the Dirichlet Laplacian defined on $\Omega$ and  by $\phi_k$  a set of orthonormal associated eigenfunctions, with $\phi_1>0$.
\begin{hyp}\label{H0} 
Assume $h \in L^q$, $q>N$ if $N\geq 2$ and $q=2$ if $N=1$. 
\end{hyp} 
\begin{hyp} \label{h1p0}
Assume  $h^1:=\int h \phi_1 >0$.
\end{hyp}
Writing 
\begin{equation} 
\label{3}h= h^1 \phi_1 + h^{\bot} 
\end{equation}
where $\int_{\Omega} h^{\perp} \phi_1 =0$  one has: 

\begin{lem} \label{L21} 
\cite{FHeT2015} We assume $ \lambda_1< \sigma \leq \Lambda <\lambda_2$ and $h \in L^q, \, q>N\geq 2$.
We suppose that there exists a  constant $C_1$ depending only on $\Omega, q,$ and $\Lambda$  such that $z$ satisfying 
(\ref{1}) is such that
\begin{equation} \label{4} \|z\|_{L^2} \leq C_1 \|h\|_{L^2}.\end{equation}
Then there exist  constants $C_2$ and $C_3$, depending only on $\Omega, q$ and $\Lambda$ such that
\begin{equation} \label{5} \|z\|_{{\cal C}^1} \leq C_2 \|h\|_{L^q} \; {\rm and  } \; \| z \|_{L^q} \leq \, C_3 \|h\|_{L^q}. \end{equation}
\end{lem}

\begin{rem} 
The same result holds for $\Lambda < \sigma < \lambda_1$ where $\Lambda$ is any given constant $< \lambda_1$, with the same proof.
\end{rem}

\begin{rem}
Inequality  (\ref{4}) cannot hold, for all $\lambda_1< \sigma \leq \Lambda$, unless $h$ is orthogonal to $\phi_1$.
\end{rem}

\begin{thm} \cite{FHeT2015}: $\,$ 
Assume Hypotheses \ref{H0} and \ref{h1p0};  fix $\Lambda$ such that  $\lambda_1<\sigma \leq  \Lambda < \lambda_2$. There exists a constant $K$ depending only on $\Omega$, $\Lambda$ and $q$ such that,
 for $\lambda_1<\sigma< \lambda_1 + \delta(h)$ with
\begin{equation} \label{ E}\delta(h)=\frac{K h^1}{\|h^{\bot}\|_{L^q}}, \end{equation}
the solution $z$ to (\ref{1}) satisfies the antimaximum principle,  that is
\begin{equation} \label{AMP }z<0 \;  {\rm in }\; \Omega; \;\; \partial z/ \partial \nu >0 \; {\rm on } \; \partial \Omega,  \end{equation}
where $\partial / \partial \nu$ denotes the outward normal derivative. 
\end{thm}

\subsection{Other  remarks  for one equation}
 Consider again Equation (\ref{1}). 
For $\sigma \neq  \lambda_k$,   $z$ solution to (\ref{1}) is
\begin{equation} \label{eqz} z\, =\, z^1 \phi_1 + z^{\perp} \,=\, \frac{h^1}{\lambda_1 - \sigma}  \phi_1 \, + \, z^{\bot},
\end{equation}
with $z^{\bot}$ satisfying
\begin{equation} \label{eqzperp}  -\Delta z^{\bot} \, =\, \sigma  z^{\bot} + \,h^{\bot} \;\, {\rm in }\; \Omega \; ; \;\;
 z^{\bot}=0 \; {\rm on } \; \partial \Omega. 
\end{equation}
In the next section, our proofs will use the  following  result. 
\begin{lem}\label{zbot} 
We assume Hypothesis \ref{H0} and  $\sigma < \lambda_1$.  
Then $z^{\perp}$ (and its first derivatives)  is bounded: 
There exits a positive constant $C_0$, independent of $\sigma$ such that
\begin{equation} \label{C0} 
 \|z^{\perp}\|_{{\cal C}^{1} } \leq C_0 \|h\|_{L^q}. 
\end{equation}
Moreover, if $\sigma < \Lambda < \lambda_1$, where $\Lambda$ is some given constant $<\lambda_1$, $z$ is bounded and there exits a positive constant $C'_0$, independent of $\sigma$ such that
\begin{equation} \label{C'0} 
 \|z\|_{{\cal C}^{1} } \leq C'_0 \|h\|_{L^q}. 
\end{equation}
\end{lem}
{\bf Proof:$\,$} 
This is a simple consequence of the 
variational characterization of $\lambda_2$:
$$ \lambda_2 \int_{\Omega} |z^{\bot}|^2 \, \leq \,  \int_{\Omega} |\nabla z^{\bot}|^2 \, = 
\, \sigma \int_{\Omega} |z^{\bot}|^2 \, +\, \int_{\Omega} z^{\bot} h^{\bot} \, \leq  \lambda_1 \int_{\Omega} |z^{\bot}|^2 \, +\, \int_{\Omega} z^{\bot} h^{\bot}.$$ 
By Cauchy-Schwarz we deduce 

\begin{equation} \label{bdd}
\| z^{\bot}\|_{L^2} \, \leq \, \frac{1}{\lambda_2-\lambda_1} \| h^{\bot}\|_{L^2}. 
 \end{equation}
This does not depend on $\sigma <\lambda_1$. 
\par \noindent
Then one can deduce  (\ref{C0}), that is $z^{\perp}$ (and its derivatives) is bounded. 
This can be found $e.g.$ in 
  \cite{Br} (for $\sigma < \lambda_1$ and $ \lambda_1 - \sigma$ small enough) or it can be derived exactly as in \cite{FHeT2015}
(where the case $\sigma > \lambda_1$  and $\sigma - \lambda_1$ small enough is considered). 
\par \noindent
Finally we write $z = z_1 \phi_1 + z^{\perp} $  and deduce (\ref{C'0}).
\begin{rem}\label{ubd}  Note that in (\ref{eqz}), since $h^1>0$, $ \frac{h^1}{\lambda_1 - \sigma} 
\rightarrow + \infty$ as $\sigma \rightarrow \lambda_1$,  $\sigma <  \lambda_1$.
\end{rem}

\section{Results for a  $n \times n$ system:} 
We consider  now a $n\times n$ (eventually non-cooperative)  system  defined on $\Omega $ a smooth bounded domain in $\RR^N$:
$$-\Delta U = AU + \mu U + F \; {\rm in }  \, \Omega \, ,  \; U=0 \; {\rm on}  \; \partial \Omega, \eqno(S)$$
where $F$ is a column vector with components $f_i$, $1 \leq i \leq n$. 
Matrix  $A$ is not necessarily cooperative, that means that its terms outside the diagonal are not necessarily positive. 
First we introduce some notations concerning matrices. Then, with these notations we can state our results and prove them. 

\subsection{The matrix of the system and  and the eigenvalues}

\begin{hyp}\label{HA} $\, $ $A$ is a $n\times n$  matrix which has constant coefficients and  has only real eigenvalues. Moreover, the largest one which is denoted by $\xi_1$ is positive and algebrically and geometrically simple.
The associated eigenvectors $X_1$ has only non zero components.  
\end{hyp}
\vskip0.2cm \noindent
Of course some of the other eigenvalues can be equal. Therefore we
write them  in decreasing order
\begin{equation} \label{order} \xi_1>\xi_2 \geq \ldots \geq \xi_n. \end{equation}
\par  \noindent
The eigenvalues of $A=(a_{ij})_{1\leq i,j\leq n}$, denoted , $\xi_1$, $\xi_2$,..., $\xi_n$ , are the roots of the associated characteristic polynomial
\begin{equation}\label{polynomial}
p_A(\xi)=det(\xi I_n-A)=\prod(\xi-\xi_k),
\end{equation}
where $I_n$ is the $n\times n$ identity matrix.
\begin{rem} By above, 
$\xi> \xi_1\Rightarrow p_A(\xi)>0$.
\end{rem}
\par \noindent
Denote by $X_1$ ... $X_n$ the  eigenvectors associated respectively to eigenvalue $\xi_1, ..., \xi_n$. 
\par \noindent
{\bf Jordan decomposition}\label{J} $\;$  Matrix A can be expressed as $A=PJP^{-1}$, where $P=(p_{ij})$ is the change of basis matrix of $A$ and $J$ is the Jordan canonical form (lower triangular matrix) associated with $A$. The diagonal entries of $J$ are the ordered eigenvalues of $A$ and  $p_A(\xi)=p_J(\xi)$.
\par \noindent
{\bf Notation} : $\,$ In the following, set 
\begin{equation}\label{tilde}  U=P \tilde U \, \Leftrightarrow \,  \tilde U= P^{-1}U, \; F=P \widetilde F \, \Leftrightarrow \,  
\tilde F = P^{-1}F.\end{equation}
Here $\widetilde U$ and $\widetilde F$ are column vectors with components 
$\widetilde u_i$ and $\widetilde f_i$.
\vskip0.2cm \noindent
{\bf Eigenvalues of the system}: $\;$  
$\mu$ is an eigenvalue of the system if there exists a non zero solution $U$ to  
$$-\Delta U = AU + \mu U  \; {\rm in }  \, \Omega \, ,  \; U=0 \; {\rm on}  \; \partial \Omega . \eqno(S_0)$$
We also say that 
$\mu$ is a "principal  eigenvalue"  of System $(S)$ if it is an eigenvalue with 
components of the associated eigenvector which does not change sign. (Note that the components do not change sign but are not necessarily positive as claimed in \cite{FHeT2015}). 
\par \noindent
Then $\phi_j X_k$ is an eigenvector associated to eigenvalue 
\begin{equation} \label{mu} \mu_{jk} = \lambda_j - \xi_k.\end{equation}

\subsection{Results  for $|\mu -  \mu_{11}| \rightarrow 0$}
\par \noindent
 We study here the sign of the component of $U$ as  $\mu \rightarrow \mu_{11} = \lambda_1- \xi_1$. 
\par \noindent
For this purpose we use  the methods in  \cite{Le} or \cite{AFLW} combined with 
\cite{FHeT2015}.
Note that  by (\ref{order}),  $\mu_{11} < \mu_{1k} = \lambda_1 - \xi_k, $ for all 
$2 \leq k \leq n$. 
\vskip0.2cm \noindent
\begin{hyp} \label{HF1}$\, $
$F$ is  with components 
$f_i \in   L^q$, $q>N>2$, $q=2$ if $N=1$,  $1\leq i \leq n$; moreover we assume that  the first   component  $\tilde f_1$  of $\tilde F =P^{-1}F$ is $\geq 0$, $\not \equiv 0$.
\end{hyp}
\begin{thm} \label{PMS} $\,$ Assume Hypothesis \ref{HA} and \ref{HF1}. Assume also  $\mu < \mu_{11}$ .   Then, there exists $\delta >0$ independant of $\mu$, such that for $ \mu_{11} - \delta < \mu <  \mu_{11}$,   the components $u_i$ of the solution $U$ have the sign of $ p_{i1}$ and  the outside normal derivatives  $\frac{\partial u_i}{\partial \nu}$ have the sign of 
 $- p_{i1}$. \end{thm}

\begin{thm} \label{AMPS} $\,$ Assume Hypothesis \ref{HA} and \ref{HF1} are satisfied; then, there exists $\delta >0$ independant of $\mu$ 
such that for   $\mu_{11} < \mu < \mu_{11} + \delta$ the components $u_i$ of the solution $U$ have the sign of $- p_{i1}$ and their outgoing normal derivatives have opposite sign.\end{thm}

\begin{rem}
The results of Theorems \ref{PMS} and \ref{AMPS} are still valid if we assume only $\int_{\Omega} \tilde f_1 \phi_1 >0$ instead of 
$ \tilde f_1 \geq 0$ $\not \equiv 0$. 
\end{rem}

\subsection{Proofs }  
We start with the proof of Theorem \ref{PMS} where $\mu < \mu_{11}$; assume Hypotheses \ref{HA}  and \ref{HF1}.
\subsubsection{Step 1: An equivalent system}
We follow \cite{Le} or \cite{AFLW}. 
As above set $U=P\tilde U$ and $F=P \tilde F$. 
\vskip0.15cm \noindent
Starting from $$-\Delta U = AU+\mu U + F,$$  multiplying by $P^{-1}$, we obtain
$$- \Delta \tilde U = J \tilde U + \mu \tilde U + \tilde F.$$
Note that everywhere we have  the homogeneous Dirichlet boundary conditions, but we do not write them for simplicity. 
\par \noindent 
The Jordan matrix $J$ has $p$ Jordan blocks $J_i$ ($1 \leq i \leq p \leq n$)
which are $k_i \times k_i$ matrices  of the form 
$$J_i=\left(\begin{array}{cccc}
 \xi_i&0&\dots&0\\
1& \xi_i& 0& \dots   \\ 
\ddots&\ddots & \vdots\\
0&\dots1&\xi_i&0 \\
0 &\dots  & 1&\xi_i
 \end{array}
 \right).$$
By Hypothesis \ref{HA}, the first block is $1\times 1$ :  $J_1=(\xi_1).$
Hence  we obtain the first equation
\begin{equation}\label{eq:9}
-\Delta \widetilde{u_1}=\xi_1\widetilde{u_1} + \mu \widetilde{u_1}  +\tilde f_1. 
\end{equation}
 
\par \noindent
Since $\tilde {f}_1\geq 0, \, \not \equiv 0$, $ \xi_1 + \mu <\lambda_1 $ and by Hypothesis \ref{HF1}, 
$\tilde{f}_1\in L^2$, we have the maximum principle and 
\begin{equation}\label{tildeu1}\widetilde {u_1}>0 \; on \; \Omega.\;\; \frac{\widetilde {u_1}}{\partial \nu}|_{\partial \Omega} <0. \end{equation}
\vskip0.15cm \noindent
Then we consider the second Jordan blocks $J_2$ which is   a $k_2\times k_2$ matrix with first line 
$$\xi_2, \,0,\,0,...$$ 
The first equation of this second block is
$$-\Delta\tilde  u_2 = \xi_2 \tilde  u_2 + \mu \tilde  u_2 + \tilde f_2.$$
Since  $\mu < \mu_{11} =\lambda_1 - \xi_1  < \lambda_1 - \xi_2 \leq \lambda_1 - \xi_k$, $k\geq 2$.  
Hence, by Lemma \ref{zbot},   $\widetilde {u_2}$ stays bounded as $\mu \rightarrow \mu_{11}$.
and this holds for all the  $\widetilde u_k$, $k>1$. 
By induction $\widetilde {u_k}$ is bounded for all $k$. 

\subsubsection{Step 2: End of the proof of Theorem \ref{PMS} }$\,$
Now we go back to the functions $u_i$: 
$U=P \tilde{U}=(u_i)$ implies that for each $u_i,\, 1\leq i\leq n$, we have

\begin{equation}\label{eq:15}
u_i=p_{i1} \widetilde{u_1} + \sum_{j=2}^n p_{ij}\widetilde {u_j}.
\end{equation}
The last term in (\ref{eq:15}) stays bounded according to Lemma \ref{zbot};  indeed
$\sum_{j=2}^n p_{ij}\widetilde {u_j}$ is bounded by a constant which does not depend on $\mu$. 
\par \noindent
By Remark \ref{ubd}, $\widetilde {u_1} \rightarrow + \infty$  as $\mu \rightarrow \lambda_1 - \xi_1$.  Hence, 
 each $u_i$ has the same sign than $p_{i1}$ (the first coefficient of the $i-th$ line in matrix $P$ which is also the $i$-th coefficient of the first eigenvector $X_1$) for $\lambda_1 - \xi_1 - \mu >0 $ small enough.  Analogously, $\frac{\partial u_i}{\partial \nu}$ behaves as $p_{i1}  \frac{\partial {\tilde u_i}}{\partial \nu}$  which has the sign of $-p_{i1} $.
\par \noindent
It is noticeable that only $\widetilde u_1$ plays a role!!
$\bullet$

\subsection{Proof of Theorem  \ref{AMPS} ( $\mu >\mu_{11}$)}
Now  $\mu_{11}<\mu < \mu_{11}+ \epsilon$ where $\epsilon \leq min\{ \xi_1 - \xi_2, \lambda_2 - \lambda_1 \}$ and $f_i \in  L^q, q>N$.   
We proceed as above but deduce immediately that for $\mu - \mu_{11}$ small enough  ($\mu - \mu_{11}  < \delta_1:= \delta(\widetilde f_1)< \frac{K\widetilde f_1^1}{\|f_1^{\perp}\|_{L^q}}$) defined in \cite{FHeT2015}, Theorem 1), 
$\widetilde u_1<0$ by the antimaximum principle. 
From now on choose 
\begin{equation}\label{delta}   \mu - \mu_{11} < \delta, \,  with \; 
\delta< min\{\epsilon, \delta_1\}.\end{equation}
For the other equations, by Lemma \ref{L21}, $\widetilde u_k >0$ is bounded  as above. 
\par \noindent
We consider now $U$. We notice that $F=P\widetilde F$ which can also be written $f_i = \sum_{k=1}^n p_{ik} \widetilde f_k$ implies
$f_i^{\perp} = \sum_{k=1}^n p_{ik} \widetilde f_k ^{\perp}$. 
With the same argument as above, 
the components $u_i$ of the  solution $U$ have the sign of $ -p_{i1}$ for $\mu - \mu_{11} $ sufficiently small  ($\mu-  \mu_{11} < \delta $). The normal derivatives of the $u_i$ are of opposite sign. 
{$\bullet$}

\section{Annex: The $2 \times 2$ non-cooperative system} We apply  now our results  to  the $2 \times 2$ system, considered in \cite{FHeT2015}. 
Consider   the $2 \times 2$ non-cooperative system depending on  a real parameter $\mu$
$$-\Delta U = AU + \mu U + F \; {\rm in }  \, \Omega \, ,  \; U=0 \; {\rm on}  \; \partial \Omega ,  \eqno(S)$$
which can also  be  written as 
$$-\Delta u \, =\, au\,+\, bv \, +\, \mu u \,  +\, f\;\, {\rm in }\; \Omega ,   \eqno(S_1)$$ 
$$-\Delta v \, =\,  cu \, + \, dv \, +\, \mu v \, +\, g \;\, {\rm in }\; \Omega,   \eqno(S_2)$$
$$ u=v=0 \; {\rm on } \; \partial \Omega.  \eqno(S_3)$$ 

\par \noindent 
\begin{hyp} \label{H1}  
Assume $b > 0 \,, c < 0,  \, $ and  
$D : = (a-d)^2 + 4bc >0. $
\end{hyp}
\vskip0.2cm \noindent
Here System $(S)$ has (at least) two principal eigenvalues $\mu^-_1 $ and $ \mu_1^+$ where 
\begin{equation} \label{mu} \mu_1^- \,: =\,  \lambda_1 - \xi_1 \,< \,  
  \mu_1^+ \,:=\, \lambda_1 - \xi_2 , \end{equation}
where $\xi_1 $ and $ \xi_2.$ are  the eigenvalues of Matrix $A$ and we choose 
$\xi_1 >\, \xi_2$.
\vskip0.2cm \noindent

\par \noindent
The main theorems in \cite{FHeT2015} are:
\begin{thm}\label{thm2}  (\cite{FHeT2015}) $\;$  Assume Hypothesis \ref{H1},  $\mu_1^- \, <  \,\mu  \,<  \,\mu_1^+$  and $  d < a  $. Assume also  
$$ f\geq 0, \, g\geq 0, \, f, g \not \equiv 0 , f,g \in L^q, \,q>N \, if \, N\geq 2 \, ;  \; q=2 \, if \, N=1.  $$
Then there exists $\delta >0$, independent of $\mu$, such that 
$ \mu < \mu_1^- + \delta $ implies 
$$ u<0, \, v>0 \; in \; \Omega; \; \frac{\partial u}{\partial \nu} >0, \, \frac{\partial v}{\partial \nu}<0 \; on \; \partial \Omega.$$ \end{thm}
\begin{thm}\label{thm4} (\cite{FHeT2015})
 $\;$  Assume Hypothesis  \ref{H1},  $\mu_1^- \, <  \,\mu  \,<  \,\mu_1^+$  and $ a< d   $. Assume also  
$$ f\leq 0, \, g\geq 0, \, f, g \not \equiv 0 \,, f,g \in L^q, \,q>N \, if \, N\geq 2 \, ;  \; q=2 \, if \, N=1.$$
Then there exists $\delta>0$, independent of $\mu$, such that i
$ \mu < \mu_1^-  + \delta$ implies 
$$
u<0, \, v<0 \; in \; \Omega; \; \frac{\partial u}{\partial \nu} >0, \, \frac{\partial v}{\partial \nu}>0 \; on \; \partial \Omega.$$ 
\end{thm}
\begin{thm} \label{thm6} (\cite{FHeT2015}) $\;$ 
Assume Hypothesis \ref{H1}  and $  a < d  $. Assume also  that the parameter $\mu$ satisfies:
 $\mu  \,<  \,\mu_1^-$ , and  
$$ f\geq 0, \,  g \geq 0, \, f, g \not \equiv 0, \, f,g \in L^2.$$
Assume also  $t ^* g -f \geq 0, \,  t ^* g -f \not \equiv 0$  with $$t^*= \frac{d-a+\sqrt D}{-2c}.$$ 
Then  
$$ u>0, \, v>0 \; in \; \Omega; \; \frac{\partial u}{\partial \nu} <0, \, \frac{\partial v}{\partial \nu}<0 \; on \; \partial \Omega.$$ 
\end{thm}  
The matrix $A$ is 
$$A=\left(\begin{array}{cccc}
a&b\\c&d
 \end{array}
 \right),$$
with  eigenvalues  $\xi_2= \frac{a+d - \sqrt D}{2} < \xi_1= \frac{a+d + \sqrt D}{2}$ where $D=(a-d)^2+4bc>0$. The eigenvectors are 
$$X_k= \left(\begin{array}{cccc}
b \\  \xi_ k  -a 
 \end{array}
 \right), \;\;  P\, =\,  \left(\begin{array}{cccc}
b & b\\ \xi_1 - a &  \xi_2 - a
 \end{array}
 \right). $$
 Note that the characteristic polynomial is ${\cal P} (s)= (a-s)(d-s) - bc$.
 Since ${\cal P}(a)={\cal P}(d)=-bc>0$, $a$ and $d$ are outside $[\xi_2, \xi_1]$. 
 \par \noindent 
For $d>a$ both $p_{i1}>0$ and for $d<a$ $p_{11}>0, \,  p_{21}<0$. 
$$P^{-1} = \frac{1}{b(\xi_1 - \xi_2)} \left(\begin{array}{cccc}
a-\xi_2&b\\\xi_1 - a & -b
 \end{array}
 \right). $$
\begin{equation} \label{cond} \tilde f_1 = \frac{1}{b(\xi_1 - \xi_2)}[ (a-\xi_2)f + bg]. \end{equation} 
\par \noindent
In Theorem 2 of \cite{FHeT2015} $d<a$, $f, g \geq 0$ so that $\tilde f_1>0$ and $u$ has the sign of $-p_{11}=-b <0$; $v$ has the sign of $-p_{21}=a-\xi_1 >0$. 
\par \noindent  
In Theorem 3 of  \cite{FHeT2015} $d>a$, $f\leq 0$ and $g \geq 0$ implies $\tilde f_1 >0$. So that $u$ has the sign of $- p_{11} = - b <0$; $v$ has the sign of $- p_{12} = a-\xi_2 <0$. 
\par \noindent  
Finally the hypothesis    $\tilde f_1 \geq 0$  is sufficient for having the sign of the solutions and the maximum principle holds (all $u_i>0$)  iff $p_{i1}>0$.
\par \noindent 
Our results can conclude for other cases; $e.g$, as in Theorem 2,  $d<a$, $f \geq 0$, but now $g<0$ with   $\tilde f_1 =  \frac{1}{b(\xi_1 - \xi_2)}[(a-\xi_2)f + bg ]>0 $. 
\par \noindent
Analogously, in Theorem 4,  $f, g \geq 0$ and   $\tilde f_1>0$ implies for having  $u,v>0$  that necessarily $\xi_2 - a >0$ so that $a<d$. But again we can conclude for the sign in other cases  ($e.g.$ $a>d$)  if only  $\tilde f_1>0$, ( which is precisely  the added condition in Theorem 4).
{$\bullet$}
\vskip0.3cm \noindent

\end{document}